\documentclass[twocolumn]{autart} 

\usepackage{graphicx}          
\usepackage{colortbl}
\usepackage{natbib}
\newtheorem{theorem}{Theorem}{}
\newtheorem{lemma}{Lemma}{}
\newtheorem{definition}{Definition}{}
{}

\newtheorem{remark}{Remark}{}
\newtheorem{assumption}{Assumption}{}

\newcommand{\bol}[1]{\boldsymbol{#1}}

\usepackage{soul,xcolor}
\usepackage{amsmath,amssymb,amsfonts}
\usepackage{mathtools}
\colorlet{shadecolor}{yellow}
\usepackage{mathrsfs}
\usepackage{array}
\usepackage{mdwtab}
\usepackage{multirow}

\date{}

\begin{document}

\begin{frontmatter}

\title{\textcolor{black}{Practical prescribed-time prescribed performance control with asymptotic convergence -- A vanishing $\bol{\sigma}$-modification approach\thanksref{footnoteinfo}}}

\thanks[footnoteinfo]{This work was supported by the Science, Technology, and Innovation Commission of Shenzhen Municipality, China under Grant No. ZDSYS20220330161800001, and the National Natural Science Foundation of China under Grant 62350055, W2433161, 62403232, and Postdoctoral Fellowship Program of CPSF under Grant No. BX20240149.\\
Corresponding author: He Kong}

\author[a]{Mehdi Golestani}\ead{golestani@sustech.edu.cn},    
\author[b]{Yongduan Song}\ead{ydsong@cqu.edu.cn},              
\author[a]{Weizhen Liu}\ead{liuwz@sustech.edu.cn}, 
\author[a,c]{Guangren Duan}\ead{duangr@sustech.edu.cn}, 
\author[a]{He Kong}\ead{kongh@sustech.edu.cn} 

\address[a]{Shenzhen Key Laboratory of Control Theory and Intelligent Systems at the Southern University of Science and Technology (SUSTech), Shenzhen, 518055, China}               
\address[b]{School of Automation, Chongqing University, Chongqing, 400044, China} 
\address[c]{Center for Control Theory and Guidance Technology, Harbin Institute of Technology, Harbin, China}    
          

\begin{abstract}                          
 In this paper, we present a method capable of ensuring practical prescribed-time control with guaranteed performance for a class of nonlinear systems in the presence of time-varying parametric and dynamic uncertainties, and uncertain control coefficients. Our design consists of two key steps.  First, we construct a performance-rate function that freezes at and after a user-specified time ${T}$, playing a crucial role in achieving desired precision within prescribed time ${T}$ and dealing with unmodeled dynamics. Next, based on this function and a ${\sigma}$–modification strategy in which the leakage term starts to vanish at ${t \geq T}$, we develop an adaptive dynamic surface control framework to reduce control complexity, deal with uncertainties, ensure prescribed performance, practical prescribed-time convergence to a specific region, and ultimately achieve asymptotic convergence. The effectiveness of the proposed control method is validated through numerical simulations.
\end{abstract}

\end{frontmatter}

\section{Introduction}
The problem of adaptive control for uncertain nonlinear systems has received significant attention in the control community. In the existing literature, backstepping control has been of interest due to its ability to deal with the case when the matching condition is not satisfied (\cite{krstic1995nonlinear, karagiannis2008nonlinear,cai2017adaptive,zhai2018output}). While most of the existing works consider nonlinear systems with parametric uncertainties (\cite{chen2021adaptive}), real-world systems might also be subject to dynamic uncertainties that can compromise the performance or even system stability (\cite{rohrs1985robustness}). The problem becomes more complicated if uncertain control coefficients are considered.

In particular, based on the backstepping control framework, substantial advancements have been achieved for nonlinear systems with parametric and/or dynamic uncertainties. However, for high-order nonlinear systems, the control complexity increases due to repeated differentiation of virtual control inputs throughout each step of the design process. Among the potential solutions to this issue is the dynamic surface control (DSC) framework, as initially introduced by \cite{swaroop2000dynamic}. Based on this concept, various DSC schemes for nonlinear systems with unknown time-varying disturbances (\cite{von2019dynamic}), parametric uncertainties (\cite{zhao2019removing}), and unmodeled dynamics (\cite{zhang2017adaptive_SMC}) have been proposed. \textcolor{black}{While \cite{zhang2017adaptive} proposed a DSC scheme considering both dynamic uncertainties and uncertain control coefficients, it can lead to an algebraic loop dilemma. Although \cite{hong_cheng_neuroadaptive_2023} addressed this issue, their method assumes a known lower bound for the control coefficient and does not ensure asymptotic convergence.}

Moreover, in the aforementioned studies, finite-time convergence of the state trajectories cannot be guaranteed. A recent advancement in the field of finite-time control has been the introduction of the concept of prescribed-time control (PsTC) by  \cite{song2017time}. Under this framework, the convergence time is established independently of initial conditions or controller parameters. Further extensions of the PsTC have been considered in the literature (see references \cite{song2017time, wang2018leader, holloway2019prescribed, holloway2019prescribed_Automatica, krishnamurthy2020dynamic,  zhou2021prescribedTAC, abel2023prescribed, zhou2023prescribed, ding2023strong, aldana2023inherent, sun2025adaptive}). Due to the use of an infinite gain, the control signal in PsTC schemes tends to infinity at the user-defined terminal time $T$. Moreover, they are only valid within the prescribed time interval $[0, T)$ and do not consider dynamic uncertainty and unknown control coefficients. 
\textcolor{black}{While \cite{krishnamurthy2023prescribed} addressed prescribed-time control problem for uncertain nonlinear systems with unknown input gain and unmodeled dynamics, their method cannot be extended to an infinite-time horizon.}


Recently, several practical PsTC frameworks with bounded control input have been proposed in the sense that the state trajectories converge to a small region and stay thereafter. In \cite{cao2022practical}, a performance function is used to guarantee a given-time given-precision (GTGP), and \cite{luo2024practical} employs a bounded time-varying scaling in a backstepping design to the same end; however, both ignore nonvanishing unmodeled dynamics and do not ensure asymptotic convergence. \textcolor{black}{\cite{shi2022prescribed} achieve global GTGP with  asymptotic tracking for strict-feedback systems, but their approach is highly complicated, requires a priori knowledge of the radius of the compact set containing the unknown parameters, and does not consider unmodeled dynamics.} To reduce the control complexity, an adaptive prescribed-time DSC algorithm for strict-feedback systems with unknown constant parameters has been proposed in \cite{zuo2023adaptive}, but it is tailored to state stabilization under vanishing uncertainties and does not handle time-varying parametric/dynamic uncertainties nor guarantee prescribed performance.

To deal with uncertain parameters, an alternative is to use the $\sigma$-modification strategy introduced by \cite{ioannou1984instability}. In this case, using a leakage term, high amplitude on the control input due to large initial conditions is prevented \cite{li2016adaptive}. To reduce the control complexity and deal with uncertain parameters in nonlinear systems, finite-time command filtered control and fixed-time DSC approaches using $\sigma$-modification have been proposed by \cite{li2019finite} and \cite{sui2023novel}, respectively. Although these two control schemes are effective, they only guarantee the ultimate boundedness of the states.  In \cite{li2022command}, a command filter tracking control method with guaranteed asymptotic convergence for nonlinear systems was proposed, utilizing a new leakage term. \textcolor{black}{However, the methods in \cite{li2016adaptive, li2022command} cannot ensure practical prescribed-time convergence and achieving such convergence is non-trivial using these approaches. Additionally, the former methods do not consider unmodeled dynamics and uncertain time-varying control coefficients, and cannot guarantee prescribed performance.}

In addition to system stability and convergence rate, providing prescribed performance is important in nonlinear systems (\cite{bikas2021tracking, lv2022prescribed, bikas2023prescribed, fotiadis2023input}). To ensure prescribed performance, \cite{bechlioulis2008robust} introduced the concept of prescribed performance control (PPC). Considering uncertain MIMO systems, a fault-tolerant control framework has been proposed by \cite{zhao2023unifying} such that both symmetric and asymmetric performance specifications are guaranteed. Based on neural networks and the DSC technique, an adaptive control scheme with full-state constraints for pure-feedback nonlinear systems under dynamic uncertainty has been proposed in \cite{hong_cheng_neuroadaptive_2023} where the algebraic loop problem is overcome. A neuroadaptive DSC-based framework for nonlinear systems under dynamic uncertainties and guaranteeing full state constraints has been designed by \cite{zhang2017adaptive}. The performance function under the existing PPC schemes has a nonzero slop at beginning which can result in a large initial control effort.

In this paper, we propose an adaptive practical prescribed-time DSC framework for nonlinear systems with uncertain control coefficients, time-varying parameters, and unmodeled dynamics. The existing methods are complicated due to repeated differentiation of virtual control laws and cannot guarantee practical prescribed-time prescribed performance with asymptotic convergence in the presence of uncertain control coefficients and nonvanishing unmodeled dynamics. Unlike conventional DSC schemes, our approach employs a  prescribed-time nonlinear filter that effectively handles these challenges, providing both prescribed performance and asymptotic convergence. The main contributions of this study are outlined as follows.

\begin{enumerate}
    \item \textcolor{black}{Based on a bounded performance-rate function, a DSC control framework containing a novel prescribed-time nonlinear filter is proposed which guarantees that the tracking errors and filter errors all enter an arbitrarily small region around the origin within a user-defined time ${T}$. }
    \item \textcolor{black}{A vanishing $\sigma$--modification strategy is incorporated into the prescribed-time DSC control framework such that practical prescribed-time stability as well as asymptotic convergence of the error trajectories are achieved even when uncertain control coefficient and nonvannishing unmodeled dynamics are considered.}
    \item \textcolor{black}{By employing a barrier Lyapunov function (BLF) approach, a novel control scheme is developed that not only guarantees prescribed performance for output tracking but also achieves practical prescribed-time convergence in the presence of nonvanishing unmodeled dynamics, without requiring strong stability assumptions on these dynamics. In addition, the use of the proposed performance function effectively reduces the initial control effort.
}
\end{enumerate}

The rest of this paper is organized as follows. The problem is formulated in \ref{se:formulation}. Section \ref{se:main} contains the main results where an adaptive practical prescribed-time prescribed performance tracking control with asymptotic convergence is proposed. Simulation results are provided in Section \ref{se:simulation}. Section \ref{se:conclusion} concludes the paper.


\section{Problem Formulation} \label{se:formulation}


\subsection{System dynamics}
Consider the following uncertain nonlinear system:
\begin{align}\label{eq:sys}
        &\dot{\xi} = q(\xi, x, t), \nonumber \\
        &\dot{x}_i = g_i(t) x_{i+1} +\theta_i^\top(t)\phi_i(\bar{x}_i)+\Delta_i(t,x,\xi), 1\leq i\leq n-1, \nonumber \\
        &\dot{x}_n = g_n(t) u+ \theta_n^\top(t)\phi_n(\bar{x}_n)+\Delta_n(t,x,\xi), 
\end{align}
where $x=\bar{x}_n=[x_1, x_2, \ldots , x_n]^\top\in\mathbb{R}^n$ is the measured state and $\bar{x}_i=[x_1, \ldots , x_i]^\top\in\mathbb{R}^i$; $y = x_1$ is the output; $\xi\in\mathbb{R}^{n_0}$ denotes the unmodeled dynamics; $u\in\mathbb{R}$ and $y\in\mathbb{R}$ are the control input and system output, respectively; $g_i(t)$ is the uncertain time-varying control coefficients; $\theta_i(t)\in\mathbb{R}^{r_i}$ is a vector of unknown time-varying parameters; and $\phi_i\in\mathbb{R}^{r_i}$ is a known and continuous function. \textcolor{black}{For simplicity, function arguments may be omitted when there is no risk of confusion throughout this work.}

\begin{assumption}\label{ass:y_d}
    \textcolor{black}{The reference trajectory and its derivatives up to second-order are known, bounded, and piecewise continuous.}
\end{assumption}

\begin{assumption} \label{ass:g}
    The sign of $g_i(t)$ is assumed to be known and there exists constants $\underline{g}_i$ and $\bar{g}_i$ satisfying $0< \underline{g}_i\leq|g_i(t)|<\bar{g}_i$ for $1\leq i\leq n$. Without loss of generality, we assume that $g_i(t)>0$.
\end{assumption}
\begin{assumption}[\cite{jiang1998design}]\label{ass:p}
    For all $(x,\xi,t)\in\mathbb{R}^n\times\mathbb{R}^{n_0}\times\mathbb{R}_+$, the inequality $$|\Delta_i(x,\xi,t)| \leq p_i^* \psi_{i1}(\|\bar{x}_i\|) + p_i^* \psi_{i2}(\|\xi\|),$$ holds, where $p_i^*$ is an unknown positive constant, $\psi_{i1}(\cdot)$ is a known nonnegative smooth function and  $\psi_{i2}(\cdot)$ is a nondecreasing continuous function such that $\psi_{i2}(0)=0$.
\end{assumption}
\begin{assumption}[\cite{jiang1998design}]\label{ass:exp}
    The unmodeled dynamics $\xi$ is exponentially input-to-state-practically stable (exp-ISpS). For system $\dot{\xi} = q(\xi,x,t)$, there exist class $\mathcal{K}_\infty$ functions ${\varpi}_1$, ${\varpi}_2$, $\Upsilon$, constants $c > 0$, $d \geq 0$, and a Lyapunov function $V(\xi)$ such that \begin{align}\label{eq:ISpS1}
    \begin{cases}
        {\varpi}_1(\|\xi\|) \leq V(\xi) \leq {\varpi}_2(\|\xi\|) & \\
        \frac{\partial V(\xi)}{\partial \xi}q(\xi,x,t) \leq -cV(\xi) + \Upsilon(|x_1|)+d. 
    \end{cases}
    \end{align}
\end{assumption}

\begin{remark}
    \textcolor{black}{Assumption \ref{ass:y_d} is standard in DSC methods. Assumption \ref{ass:g} ensures that the sign of the control gain \( g_i(t) \) is known and \( g_i(t) \) is nonzero, which is essential for guaranteeing controllability of the system. 
    Assumption \ref{ass:p} bounds the dynamic uncertainties \( \Delta_i(t, x, \xi) \) using functions of the states and unmodeled dynamics. It allows the adaptive control law to compensate for these nonvanishing uncertainties. Assumption \ref{ass:exp} ensures that the \( \xi \)-system is exp-ISpS with respect to the input \( x \). A consequence of this property is that the \( \xi \)-system exhibits bounded-input bounded-state stability. Specifically, when both \( x_1 = 0 \) and \( d = 0 \), the \( \xi \)-system becomes uniformly exponentially stable at \( \xi = 0 \).}
\end{remark}
\begin{lemma}[\cite{jiang1998design}]\label{lem:exp}
    If $V$ is an exp-ISpS Lyapunov function for the system $\dot{\xi}=q(\xi,x,t)$, i.e., the conditions in \eqref{eq:ISpS1} hold, then, for any constant $\bar{c}\in(0, c)$, initial instant $t_0 > 0$, initial condition $\xi_0 = \xi(t_0)$, continuous function $\bar{\Upsilon}(x_1)$ such that $\bar{\Upsilon}(x_1)\geq {\Upsilon}(|x_1|)$, there exist a finite $T^0 = \max\{0, \log \big(\frac{V(\xi_0)}{r_0}\big)/(c-\bar{c})\}\geq0$, a nonnegative function $D(t_0, t)$, and a signal described by $$\dot{r} = -\bar{c}r+\bar{\Upsilon}(x_1)+d,~ r(t_0)=r_0,$$
such that $D(t_0, t) = 0$ for $t \geq t_0 + T_0$, and $V (\xi) \leq r(t) + D(t_0, t)$ with $D(t_0, t) = \max\{0, e^{-c(t-t_0)} V (\xi_0)- e^{-\bar{c}(t-t_0)}r_0\}$.
\end{lemma}
\begin{remark}\label{remark:xi}
    \textcolor{black}{From Assumption \ref{ass:exp}, we have $\|\xi\| \leq \varpi_1^{-1}(V(\xi))$. Based on Lemma \ref{lem:exp}, there exists a positive constant $D$ such that $\|\xi\| \leq \varpi_1^{-1}(r+D)$, $\forall t\geq 0$. This inequality is used in \eqref{eq:sec1} to deal with dynamic uncertainties (\cite{zhang2017adaptive}).}
\end{remark}

\begin{lemma}[\cite{jiang1998design}]\label{lem:g}
    For any $\tau>0$, there exists a smooth, odd function $g(x)$ such that $xg(x)>0$, $g(0)=0$ and $|x| \leq x g(x) + \tau, ~\forall x\in \mathbb{R}.$  
\end{lemma}
\begin{lemma}[\cite{jiang1998design}]\label{lem:f}
    For any $\tau>0$ and any continuous function $f:\mathbb{R}\to\mathbb{R}$ with $f(0)=0$, there is a nonnegative smooth function $\hat{f}$ with $\hat{f}(0)=\partial\hat{f}/\partial x(0)=0$ such that $|f(x)| \leq \hat{f}(x) + \tau, ~\forall x\in \mathbb{R}.$
\end{lemma}

\begin{lemma}[\cite{li2022command}] \label{lem:sqrt}
    For any scalar $s\in\mathbb{R}$, vectors $\Theta\in\mathbb{R}^n$ and $\Phi\in\mathbb{R}^n$ with $\|\Theta\| \leq \vartheta$, and positive integrable time-varying function $\tau$, the following inequality holds $$s \Theta^\top \Phi \leq \vartheta \frac{s^2 \Phi^\top\Phi}{\sqrt{s^2 \Phi^\top\Phi + \tau^2}} + \tau \vartheta.$$
\end{lemma}

\begin{lemma}[\cite{ren2010adaptive}]\label{lem:log}
    For any $s\in\mathbb{R}$ and positive constant $k$, the  inequality $\log\tfrac{k^2}{k^2-s^2}<\tfrac{s^2}{k^2-s^2}$ is satisfied in the set $|s|<k$.
\end{lemma}

\subsection{Control Objective} 
For the complex time-varying system \eqref{eq:sys}, the control objective is to develop an adaptive prescribed-time DSC scheme with prescribed performance and asymptotic convergence. The designed controller ensures that
\begin{enumerate}
    \item \textcolor{black}{The output tracking error}, virtual errors, and filter errors converge to an arbitrarily small region within a prescribed time and subsequently asymptotically are regulated to zero; 
    \item The output trajectory always remains within the prescribed boundary; and
    \item All signals in the controlled system are bounded.
\end{enumerate}


\section{Main results} \label{se:main}


\subsection{\textcolor{black}{Performance-Rate Function}}
\textcolor{black}{To streamline the design, we introduce a performance-rate function $\mu(t)$ that can be utilized as either performance function $\rho(t)$ or rate function $\sigma(t)$.}
\textcolor{black}{
\begin{definition}\label{def:mu}
A $\textcolor{black}{C^1}$ function $\textcolor{black}{\mu(t): [0, \infty) \to (0, \infty)}$ is called a performance-rate function if for user-defined constants $\textcolor{black}{\mu_0 > 0}$, $\textcolor{black}{\mu_T > 0}$, and a prescribed time $\textcolor{black}{T > 0}$, it satisfies the following properties:
\begin{enumerate}
\item $\textcolor{black}{\mu(0) = \mu_0}$ and $\textcolor{black}{\mu(t) = \mu_T}$ for all $\textcolor{black}{t \geq T}$.
\item It is strictly increasing on $\textcolor{black}{[0, T)}$ if $\textcolor{black}{\mu_T > \mu_0}$, and strictly decreasing if $\textcolor{black}{\mu_T < \mu_0}$.
\item $\textcolor{black}{\dot{\mu}(0) = 0}$.
\end{enumerate}
\end{definition}
\begin{remark}\label{remark:mu-example}
An example of $\textcolor{black}{\mu(t)}$ is as follows:
\begin{align}\label{eq:mu}
\textcolor{black}{
\mu(t) = \begin{cases}
\mu_T + \dfrac{(\mu_0 - \mu_T)(T-t)^2}{(\upsilon_\mu^2+1)t^2 + T^2 - 2T t}, & t \in [0, T) \\
\mu_T, & t \in [T, \infty)
\end{cases}
}
\end{align}
where $\textcolor{black}{\upsilon_\mu > 0}$ is a design parameter.
\end{remark}}

\begin{remark}
    \textcolor{black}{We define two key functions from $\mu(t)$: the performance function defined as $\rho(t) = \mu(t)$ with $\mu_0 = \rho_0$ and $\mu_T = \rho_T$, where $0 < \rho_T < \rho_0$ ensuring $\rho(t)$ is decreasing; and the rate function defined as $\sigma(t) = \mu(t)$ with $\mu_0 = 1$ and $\mu_T = \bar{\sigma}>1$, ensuring $\sigma(t)$ is increasing.}
\end{remark}

\begin{definition}\label{def:eps}
    A function $\varepsilon(t)>0$ is integrable and differentiable if $\int_0^t\varepsilon(s)ds\leq \bar{\varepsilon}<\infty$ and $\dot{\varepsilon}(t)\leq \underline{\varepsilon}<\infty$.
\end{definition}


\subsection{Controller design}

To design the adaptive DSC scheme, let us define the following coordinate transformations:
\begin{align}
    z_1 &= (x_1-\alpha_{0}^c)/\rho, ~ z_i = x_i-\alpha_{i-1}^c;~i= 2, 3, \ldots, n, \label{eq:z}\\
    \omega_j &= \alpha_j^c-\alpha_j;\hspace{0.6in} j = 1, 2, \ldots, n-1, \label{eq:w}
\end{align}
where $\alpha_{0}^c=y_d$ is the desired trajectory, $z_i$ is the virtual error, $\alpha_j$ is the virtual control law passed through a nonlinear filter to obtain $\alpha_j^c$, and $\omega_j$ is the output error of the filter.

\textbf{Step 1:} Based on Eqs. \eqref{eq:z}, \eqref{eq:w} and \eqref{eq:sys}, the dynamics of $z_1$ can be given as
\begin{align}\label{eq:z1}
    \dot{z}_1 = \frac{g_1( z_2+\omega_1+\alpha_1)+\theta_1^\top\phi_1-\dot{y}_d+\Delta_1}{\rho}-\frac{\dot{\rho}}{\rho}z_1.
\end{align}
According to Assumption \ref{ass:p}, one has
\begin{align}\label{eq:z1mu1del1}
    z_1 \tfrac{1}{\rho} \Delta_1 \leq p_1^* \tfrac{1}{\rho} |z_1|\psi_{11}(|x_1|) + p_1^* \tfrac{1}{\rho} |z_1|\psi_{12}(\|\xi\|).
\end{align}
Using Lemmas \ref{lem:g} and \ref{lem:f}, the first term in the right side of \eqref{eq:z1mu1del1} can be written as 
\begin{align}\label{eq:1_z1mu1del1}
    p_1^* \tfrac{1}{\rho} |z_1|\psi_{11}(|x_1|) \leq p \tfrac{1}{\rho} z_1\hat{\psi}_{11}(x_1) + p\tfrac{1}{\rho} {\varepsilon},
\end{align}
where $p=\max\{p_1^*, \ldots , p_n^*\}$, $\hat{\psi}_{11}(\cdot)$ is a smooth function with $\hat{\psi}_{11}(0)=0$ and $\varepsilon$ is defined in Definition \ref{def:eps}. 

\textcolor{black}{Based on Remark \ref{remark:xi}}, one obtains
\begin{align} \label{eq:sec1}
    \psi_{12}(\|\xi\|) \leq \psi_{12} \circ \varpi_1^{-1}(r+D).
\end{align}
Since $\psi_{12} \circ \varpi_1^{-1}$ is a nonnegative nondecreasing function, then the following inequality holds
\begin{align}\label{eq:sec2}
    \psi_{12} \circ \varpi_1^{-1}(r+D) \leq \psi_{12} \circ \varpi_1^{-1}(2r) + \psi_{12} \circ \varpi_1^{-1}(2D).
\end{align}
From \eqref{eq:sec1}-\eqref{eq:sec2} and using Young's inequality, we have
\begin{align*}
    p_1^* \tfrac{1}{\rho} |z_1|\psi_{12}(\|\xi\|) 
    \leq & p\tfrac{1}{\rho}|z_1|\psi_{12} \circ \varpi_1^{-1}(2r)+\tfrac{1}{4\rho^2}z_1^2 +d_0,
\end{align*}
where $d_0=(p \psi_{12} \circ \varpi_1^{-1}(2D))^2$. Using Lemmas \ref{lem:g} and \ref{lem:f}, there exists a smooth function $\hat{\psi}_{12}(z_1,r)$ which is zero at zero such that
\begin{align*}
   p\tfrac{1}{\rho}|z_1|\psi_{12} \circ \varpi_1^{-1}(2r) \leq  p\tfrac{1}{\rho} z_1\hat{\psi}_{12}(z_1,r) +p\tfrac{1}{\rho}\varepsilon.
\end{align*}
Thus, the second term in the right side of \eqref{eq:z1mu1del1} can be written as 
\begin{align}\label{eq:2_z1mu1del1}
    p_1^* \tfrac{1}{\rho} |z_1|\psi_{12}(\|\xi\|) 
    \leq & \tfrac{p}{\rho} z_1 \hat{\psi}_{12}(z_1, r)+\tfrac{p}{\rho}\varepsilon + \tfrac{z_1^2}{4\rho^2} + d_0.
\end{align}
Then, from Eqs. \eqref{eq:z1}, \eqref{eq:1_z1mu1del1} and \eqref{eq:2_z1mu1del1}, one obtains
\begin{align}\label{eq:z1_2}
    z_1\dot{z}_1 \leq  & \tfrac{1}{\rho}z_1  \big( g_1 z_2 + g_1  \omega_1+g_1\alpha_1+\theta_1^\top\phi_1-\dot{y}_d \big)  + \tfrac{1}{4\rho^2}z_1^2 \nonumber\\
    & + p\tfrac{1}{\rho}z_1(\hat{\psi}_{11}+\hat{\psi}_{12}) - \tfrac{\dot{\rho}}{\rho}z_1^2 + D_1(t_0,t),
\end{align}
where $D_1(t_0,t)=\tfrac{p}{\rho}\varepsilon + d_0$. Based on the BLF approach, we define 
\begin{align*}
    V_1 = \frac{1}{2\underline{g}_1}\log\frac{1}{1-z_1^2}+\frac{1}{2}\omega_1^2+\frac{1}{2\iota_{\vartheta_1}}\tilde{\vartheta}_1^2 +\frac{1}{2\iota_{\gamma_1}}\tilde{\gamma}_1^2,
\end{align*}
where $\tilde{*} = * - \hat{*}$, $\iota_{\vartheta_1}, \iota_{\gamma_1}>0$, $\vartheta_1$ is an unknown constant and $\gamma_1$ is the upper bound of the first virtual control. 

Taking derivative of $V_1$ and using Eq. \eqref{eq:z1_2}, one has
\begin{align} \label{eq:V1dot}
    \dot{V}_1     \leq & \kappa_1 z_1 \zeta_1 + \kappa_1 g_1z_1z_2 + \kappa_1 g_1z_1\alpha_1 +\kappa_1 \omega_1^2 -\tfrac{1}{\iota_{\vartheta_1}}\tilde{\vartheta}_1\dot{\hat{\vartheta}}_1  \nonumber \\
    &  + \kappa_1 z_1 \Big(\tfrac{{g}_1^2 z_1}{4} + \theta_1^\top\phi_1  + p(\hat{\psi}_{11}+\hat{\psi}_{12})\Big)  \nonumber\\
    &  + \omega_1 \dot{\omega}_1 -\frac{1}{\iota_{\gamma_1}} \tilde{\gamma}_1\dot{\hat{\gamma}}_1 + \frac{\lambda}{\underline{g}_1} D_1(t_0,t)
\end{align}
where $\lambda=\frac{1}{1-z_1^2}$, $\zeta_1 = \tfrac{1}{4\rho}z_1-\dot{y}_d -\dot{\rho} z_1$, and \textcolor{black}{$\kappa_1=\frac{\lambda}{\underline{g}_1\rho}$.} Let ${\Theta}_1 = \big[\theta_1^\top, p, {g}_1^2 \big]^\top$, $\Phi_1 = \big[\phi_1^\top, (\hat{\psi}_{11}+\hat{\psi}_{12}), \tfrac{ z_1}{4}\big]^\top$, and  ${\vartheta}_1 = \sup_{t\geq 0}\|{\Theta}_1(t)\|$.
According to Lemma \ref{lem:sqrt}, one has
\begin{align}
    \textcolor{black}{\kappa_1} z_1 \Theta_1^\top \Phi_1  \leq \vartheta_1 \varepsilon + \vartheta_1 \textcolor{black}{\kappa_1}z_1\varphi_{1}, \label{eq:ineq2}
\end{align}
where $\varphi_{1} = \frac{\textcolor{black}{\kappa_1}z_1\Phi_1^\top\Phi_1}{\sqrt{\textcolor{black}{\kappa_1^2}z_1^2\Phi_1^\top\Phi_1+\varepsilon^2}}$, and $\varepsilon$ is defined in Definition \ref{def:eps}.
Substituting \eqref{eq:ineq2} into \eqref{eq:V1dot}, we have
\begin{align} \label{eq:V1dot2}
    \dot{V}_1 \leq & -\tfrac{\varsigma_{z_1}\sigma_1}{\underline{g}_1} \lambda z_1^2 + \textcolor{black}{\kappa_1}z_1\bar{\alpha}_1 + \textcolor{black}{\kappa_1}g_1z_1z_2 + \textcolor{black}{\kappa_1}g_1z_1\alpha_1  \nonumber\\
    & +\textcolor{black}{\kappa_1}\omega_1^2 + \vartheta_1 \varepsilon + \tfrac{1}{\iota_{\vartheta_1}} \tilde{\vartheta}_1 \big(\iota_{\vartheta_1}\textcolor{black}{\kappa_1}z_1 \varphi_1 - \dot{\hat{\vartheta}}_1\big)  \nonumber\\
    & + \omega_1\dot{\omega}_1 -\tfrac{1}{\iota_{\gamma_1}} \tilde{\gamma}_1\dot{\hat{\gamma}}_1 + \tfrac{\lambda}{\underline{g}_1} D_1(t_0,t),
\end{align}
where $\varsigma_{z_1}>\tfrac{1}{2}$ and the term $\kappa_1 z_1\bar{\alpha}_1$ is added and subtracted where $\bar\alpha_1$ defined in Eq. \eqref{eq:alpha1}. The first virtual control law, the update law for $\hat{\vartheta}_1$, the nonlinear filter, and the update law for $\hat{\gamma}_1$ are, respectively, developed as
\begin{align}
    &\textcolor{black}{\alpha_1 = -\frac{\textcolor{black}{\kappa_1}z_1 \bar{\alpha}_1^2}{\textcolor{black}{\underline{g}_1}\sqrt{(\textcolor{black}{\kappa_1}z_1 \bar{\alpha}_1)^2+ \varepsilon^2}}};
    \bar{\alpha}_1 = \varsigma_{z_1}\sigma_1 \rho \textcolor{black}{z_1} +\textcolor{black}{\zeta_1} + \hat{\vartheta}_1\varphi_1, \label{eq:alpha1} \\
    & \dot{\alpha}_1^c = -\varsigma_{\omega_1}\sigma_1\omega_1 - \textcolor{black}{\kappa_1\omega_1 }- \frac{\hat{\gamma}_1^2 \omega_1}{\sqrt{\hat{\gamma}_1^2\omega_1^2+\varepsilon^2}}; ~~\alpha_1^c(0)=\alpha_1(0), \label{eq:filter1} \\
    &\dot{\hat{\vartheta}}_1 = \iota_{\vartheta_1} \textcolor{black}{\kappa_1}z_1 \varphi_1 - 2\iota_{\vartheta_1}\sigma_2\hat{\vartheta}_1, \label{eq:theta1}  \\
    &\dot{\hat{\gamma}}_1 = \iota_{\gamma_1}|\omega_1| -2\iota_{\gamma_1}\sigma_2\hat{\gamma}_1, \label{eq:gamma1}
\end{align}
where $ \varsigma_{\omega_1}>0$, $\sigma_1={\sigma}$ is the rate function in Definition \ref{def:mu}, $\sigma_2={\sigma}$ if $t<T$ and $\sigma_2={\sigma}\varepsilon(t-T)$ if $t\geq T$, and $\varepsilon(t)$ is given in Definition \ref{def:eps}.

\begin{figure}
  \centering
  \includegraphics[width=2.1in]{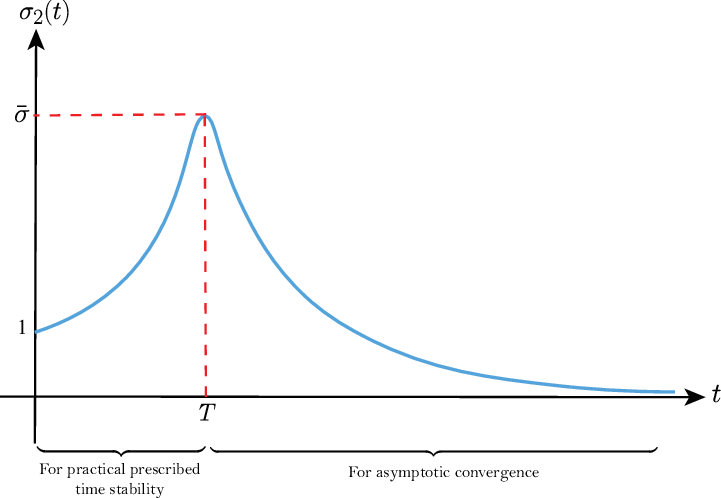}
  \caption{Schematic illustration of the time-varying gain $\sigma_2(t)$} \label{fig:sigma2}
\end{figure}

\begin{remark}
   The gain $\sigma_2(t)$ in the leakage term of update laws \eqref{eq:theta1}--\eqref{eq:gamma1} reaches $\Bar{\sigma}$ to provide practical prescribed-time stability and then starts to vanish after $t\geq T$ to guarantee asymptotic convergence. A graphical illustration of the gain $\sigma_2(t)$ is provided in Fig. \ref{fig:sigma2}. 
\end{remark}

According to Assumption \ref{ass:g}, we have $g_1(t)\geq \underline{g}_1 > 0$. From Eq. \eqref{eq:alpha1}, one has
\begin{align}\label{eq:alpha_eps}
     & \textcolor{black}{\kappa_1 g_1z_1\alpha_1  + \textcolor{black}{\kappa_1} z_1\bar{\alpha}_1 \leq   \tfrac{\big(1 - \tfrac{g_1}{\underline{g}_1} \big)\big(\textcolor{black}{\kappa_1}z_1 \bar{\alpha}_1\big)^2}{\sqrt{\big(\textcolor{black}{\kappa_1} z_1\bar{\alpha}_1\big)^2+ \varepsilon^2}} + \varepsilon \leq \varepsilon.}
 \end{align}
Then, substituting \eqref{eq:alpha1}--\eqref{eq:alpha_eps} into \eqref{eq:V1dot2}, one obtains
\begin{align} \label{eq:V1dot3}
    \dot{V}_1 \leq &  \textcolor{black}{\kappa_1}g_1z_1z_2 -\tfrac{\varsigma_{z_1}\sigma_1}{\underline{g}_1} \lambda z_1^2 - \sigma_2 (\tilde{\vartheta}_1^2-{\vartheta}_1^2) - \varsigma_{\omega_1}\sigma_1 \omega_1^2  \nonumber\\
    &  - \sigma_2 (\tilde{\gamma}_1^2-{\gamma}_1^2)  + (\vartheta_1+2) \varepsilon  + \tfrac{\lambda}{\underline{g}_1} D_1(t_0,t),
\end{align}
where $2\sigma_2 \tilde\vartheta_1\hat\vartheta_1 \le -\sigma_2\tilde\vartheta_1^2 + \sigma_2 \vartheta_1^2$ and $2\sigma_2 \tilde\gamma_1\hat\gamma_1 \le -\sigma_2\tilde\gamma_1^2 + \sigma_2 \gamma_1^2$ have been used.


\textcolor{black}{\textbf{Step \textit{i} $(2\leq i \leq n-1):$}} From \eqref{eq:z}, $z_i$ dynamics can be given as
\begin{align}\label{eq:z2dot}
    \dot{z}_i = g_i(z_{i+1}+\omega_i+\alpha_i)+ \theta_{i}^\top\phi_i + \Delta_i - \dot{\alpha}_{i-1}^c.
\end{align}
Following the same procedure in Step 1, one has
\begin{align*}
    z_i \Delta_i \leq & p z_i (\hat{\psi}_{i1} + \hat{\psi}_{i2}) + 2p{\varepsilon}   + \frac{1}{4}z_i^2 + (p \psi_{i2} \circ \varpi_1^{-1}(2D))^2.
\end{align*}
Define a Lyapunov function candidate as
\begin{align*}
    V_i = V_{i-1} +\frac{1}{2\underline{g}_i} z_i^2 + \frac{1}{2\iota_{\vartheta_i}}\tilde{\vartheta}_i^2 +\frac{1}{2}\omega_i^2+\frac{1}{2\iota_{\gamma_i}}\tilde{\gamma}_i^2,
\end{align*}
where $\iota_{\vartheta_i}, \iota_{\gamma_i}>0$. Derivative of $V_i$ along \eqref{eq:z2dot} is
\begin{align}\label{eq:V2dot}
    \dot{V}_i  \leq & \textcolor{black}{\kappa_{i-1}} g_{i-1}z_i z_{i-1} -\tfrac{\varsigma_{z_1}\sigma_1}{\underline{g}_1}  \lambda z_1^2 -\sum_{j=2}^{i-1}\tfrac{\varsigma_{z_j}\sigma_1}{\underline{g}_j} z_j^2  \nonumber\\
    & -\sum_{j=1}^{i-1}\varsigma_{\omega_j}\sigma_1 \omega_j^2 -  \sigma_2\sum_{j=1}^{i-1}\big(\tilde{\vartheta}_j^2-\vartheta_j^2 + \tilde{\gamma}_j^2- \gamma_j^2\big)  \nonumber\\
    & + \textcolor{black}{g_i\kappa_i} z_i\omega_i +\textcolor{black}{g_i\kappa_i}z_{i+1}z_i +\sum_{j=1}^{i-1}(\vartheta_j+2)\varepsilon    +\textcolor{black}{g_i\kappa_i}z_i \alpha_i  \nonumber\\
    &  + \textcolor{black}{\kappa_i}z_i\theta_i^\top\phi_i+\textcolor{black}{\kappa_i}p z_i(\hat{\psi}_{i1}+\hat{\psi}_{i2})+ \textcolor{black}{\kappa_i}z_i\zeta_i -\tfrac{1}{\iota_{\vartheta_i}}\tilde{\vartheta_i}\dot{\hat{\vartheta}}_i \nonumber\\
    & +\omega_i\dot{\omega}_i-\tfrac{1}{\iota_{\gamma_i}}\tilde{\gamma}_i\dot{\hat{\gamma}}_i+D_i(t_0,t),
\end{align}
where $\zeta_i = \frac{1}{4}z_i-\dot{\alpha}_{i-1}^c$, $\textcolor{black}{\kappa_i=\tfrac{1}{\underline{g}_i}}$, and $D_i= \frac{\lambda}{\underline{g}_1} D_1 + \sum_{j=2}^{i} \kappa_j \big(2p\varepsilon+(p\psi_{j2}\circ\varpi_1^{-1}(2D))^2)$. Let ${\Theta}_{i} = \big[g_{i-1}, \theta_i^\top, $ $ p, g_i^2\big]^\top$, $\Phi_i = \Big[\frac{\kappa_{i-1}}{\kappa_i}z_{i-1}, \phi_i^\top, \hat{\psi}_{i1}+\hat{\psi}_{i2}, \tfrac{z_i}{4}\Big]^\top$, and ${\vartheta}_i = \sup_{t\geq 0}\|{\Theta}_i(t)\|$.
Based on Lemma \ref{lem:sqrt}, one has
\begin{align}
    \textcolor{black}{\kappa_i}z_i \Theta_i^\top \Phi_i  \leq \vartheta_i \varepsilon + \vartheta_i \textcolor{black}{\kappa_i}z_i\varphi_{i}, \label{eq:ineq22}
\end{align}
where $\varphi_{i} = \frac{\textcolor{black}{\kappa_i}z_i\Phi_i^\top\Phi_i}{\sqrt{\textcolor{black}{\kappa_i^2}z_i^2\Phi_i^\top\Phi_i+\varepsilon^2}}$.

The virtual control law $\alpha_i$, the $i$th nonlinear filter, the update law for $\hat{\theta}_i$, and the update law for $\hat{\gamma}_i$ are, respectively, developed as
\begin{align}
    \textcolor{black}{\alpha_i} & \textcolor{black}{= -\frac{\textcolor{black}{\kappa_i}z_i \bar{\alpha}_i^2}{\textcolor{black}{\underline{g}_i}\sqrt{(\textcolor{black}{\kappa_i}z_i \bar{\alpha}_i)^2+ \varepsilon^2}}};~~ \textcolor{black}{\bar{\alpha}_i = \varsigma_{z_i}\sigma_1z_i +\textcolor{black}{\zeta_i} + \hat{\vartheta}_i\varphi_i, }\label{eq:alphai} \\
    \dot{\alpha}_i^c &= -\varsigma_{\omega_i}\sigma_1\omega_i-\textcolor{black}{\kappa_i \omega_i}-\frac{\hat{\gamma}_i^2 \omega_i}{\sqrt{\hat{\gamma}_i^2 \omega_i^2+\varepsilon^2}}, ~~\alpha_i^c(0)=\alpha_i(0),\label{eq:filteri} \\
    \dot{\hat{\vartheta}}_i &= \iota_{\vartheta_i}\textcolor{black}{\kappa_i} z_i \varphi_i - 2\iota_{\vartheta_i}\sigma_2\hat{\vartheta}_i, \label{eq:thetai}  \\
    \dot{\hat{\gamma}}_i &= \iota_{\gamma_i}|\omega_i|-2\iota_{\gamma_i}\sigma_2\hat{\gamma}_i, \label{eq:gammai}
\end{align}
where $\varsigma_{z_i}, \varsigma_{\omega_{i}}>0$. Considering \eqref{eq:ineq22} and substituting \eqref{eq:alphai}--\eqref{eq:gammai} into \eqref{eq:V2dot}, we have
\begin{align}\label{eq:V2dot2}
    \dot{V}_i  \leq &  \textcolor{black}{\kappa_i}g_i z_{i+1}z_i -\tfrac{\varsigma_{z_1}\sigma_1}{\underline{g}_1} \lambda z_1^2 -\sum_{j=2}^{i}\tfrac{\varsigma_{z_j}\sigma_1}{\underline{g}_i} z_j^2 +D_i(t_0,t) \nonumber\\
    & -\sum_{j=1}^{i}\varsigma_{\omega_j}\sigma_1 \omega_j^2 - \sigma_2 \sum_{j=1}^{i}\big(\tilde{\vartheta}_j^2-\vartheta_j^2+\tilde{\gamma}_j^2- \gamma_j^2\big)  \nonumber\\
    & +\sum_{j=1}^{i}(2\vartheta_j+2)\varepsilon.
\end{align}
\textbf{Step n:} From \eqref{eq:z}, $z_n$ dynamics can be given as
\begin{align*}
    \dot{z}_n = g_n u + \theta_n^\top\phi_n + \Delta_n - \dot{\alpha}_{n-1}^c.
\end{align*}
Following the same procedure, we design the control law and the update law as
\begin{align}
    &\textcolor{black}{u = -\frac{\textcolor{black}{\kappa_n}z_n \bar{\alpha}_n^2}{\textcolor{black}{\underline{g}_n}\sqrt{(\textcolor{black}{\kappa_n}z_n \bar{\alpha}_n)^2+ \varepsilon^2}};~\bar{\alpha}_n = \varsigma_{z_n}\sigma_1z_n + \textcolor{black}{\zeta_n} +\hat{\vartheta}_n\varphi_n,} \label{eq:u}  \\
    &\dot{\hat{\vartheta}}_n = \iota_{\vartheta_n} \textcolor{black}{\kappa_n} z_n \varphi_n - 2\iota_{\vartheta_n}\sigma_2\hat{\vartheta}_n, \label{eq:thetan}
\end{align}
where $\varsigma_{z_n}>0$.
Define a Lyapunov function candidate as $V_n = V_{n-1} +\frac{1}{2\underline{g}_n} z_n^2 + \frac{1}{2\iota_{\vartheta_n}}\tilde{\vartheta}_n^2$.
The, we have
\begin{align} \label{eq:Vndot}
    \dot{V}_n \leq & \textcolor{black}{\kappa_{n-1}}g_{n-1}z_n z_{n-1}  -\tfrac{\varsigma_{z_1}\sigma_1}{\underline{g}_1} \lambda z_1^2 -\sum_{i=2}^{n-1}\tfrac{\varsigma_{z_i}\sigma_1}{\underline{g}_i} z_i^2  \nonumber\\
    & -\sum_{i=1}^{n-1} \varsigma_{\omega_i}\sigma_1 \omega_i^2 - \sigma_2 \sum_{i=1}^{n-1}  \big(\tilde{\vartheta}_i^2 - {\vartheta}_i^2+\tilde{\gamma}_i^2-{\gamma}_i^2\big)  \nonumber\\
    &  +\sum_{i=1}^{n-1}(\vartheta_i+2)\varepsilon + \textcolor{black}{\kappa_n} z_n \zeta_n   +D_n + \textcolor{black}{\kappa_n}g_n z_n  u \nonumber\\
    & + \textcolor{black}{\kappa_n}z_n\Theta_n^\top\Phi_n - \frac{1}{\iota_{\vartheta_n}}\tilde{\vartheta}_n\dot{\hat{\vartheta}}_n,
\end{align}
where $\zeta_n = \frac{1}{4}z_n-\dot{\alpha}_{n-1}^c$, $\textcolor{black}{\kappa_n=\tfrac{1}{\underline{g}_n}}$, and $D_n=  D_{n-1} + \frac{1}{\underline{g}_{n}}\big(2p\varepsilon+(p\psi_{n2}\circ\varpi_1^{-1}(2D))^2\big)$. Let ${\Theta}_{n} = \big[g_{n-1}, \theta_n^\top, p, g_n^2\big]^\top$, $\Phi_n = \Big[\frac{\kappa_{n-1}}{\kappa_n}z_{n-1}, \phi_n^\top, \hat{\psi}_{n1}+\hat{\psi}_{n2}, \tfrac{z_n}{4}\Big]^\top$, and ${\vartheta}_n = \sup_{t\geq 0}\|{\Theta}_n(t)\|$.

Substituting \eqref{eq:u} and \eqref{eq:thetan} into \eqref{eq:Vndot} gives
\begin{align}\label{eq:Vndot2}
    \dot{V}_n \leq & -\tfrac{\varsigma_{z_1}\sigma_1}{\underline{g}_1} \lambda z_1^2 -\sum_{i=2}^{n}\tfrac{\varsigma_{z_i}\sigma_1}{\underline{g}_i} z_i^2 - \sum_{i=1}^{n} \sigma_2 (\tilde{\vartheta}_i^2-{\vartheta}_i^2) \nonumber\\
    & -\sum_{i=1}^{n-1} \varsigma_{\omega_i}\sigma_1 \omega_i^2 - \sum_{i=1}^{n-1}  \sigma_2 (\tilde{\gamma}_i^2- {\gamma}_i^2)   + \tilde{D},
\end{align}
where $\tilde{D} = \sum_{i=1}^{n-1}(2\vartheta_i+2)\varepsilon +(2\vartheta_n+1)\varepsilon  +D_n$.

\begin{theorem}\label{the1}
    Consider the system \eqref{eq:sys} under Assumptions \ref{ass:y_d}--\ref{ass:exp} with the control input \eqref{eq:u}, the virtual control inputs \eqref{eq:alpha1} and \eqref{eq:alphai}, the update laws for $\hat{\vartheta}_i$ in \eqref{eq:theta1}, \eqref{eq:thetai}, \eqref{eq:thetan}, the updates laws for $\hat{\gamma}_i$ in \eqref{eq:gamma1}, \eqref{eq:gammai}, and the nonlinear filters \eqref{eq:filter1}, \eqref{eq:filteri}. Given initial value $V_n(0)\leq q$, the following results are established. 1) All closed-loop system signals are bounded. 2) The output evolves within the prescribed set $\Omega_e=\{(e, \rho)~|~|e(t)|<\rho(t)\}$. 3) \textcolor{black}{The output tracking error,} virtual errors and filter errors converge to arbitrary small regions in a specified time $T$ and asymptotically converge to zero.
    \end{theorem}

\textbf{Proof.}
        1) Define compact sets $\Omega_{\rho}=\{[\rho, \dot{\rho}, \ddot{\rho}]^\top : \rho+ \dot{\rho}+ \ddot{\rho}\leq \mathcal{B}_1\}\subset \mathbb{R}^3$ and $\Omega_V=\big\{\tfrac{1}{\underline{g}_1}\log\tfrac{1}{1-z_1^2}+\sum_{i=2}^{n}\tfrac{1}{\underline{g}_i}z_i^2 +\sum_{i=1}^{n}\tfrac{1}{\iota_{\vartheta_i}}\tilde{\vartheta}_i^2 + \sum_{i=1}^{n-1}\big(\omega_i^2+\tilde{\gamma}_i^2\big)\leq 2q\big\}\subset \mathbb{R}^{4n-2}$. Therefore, $\log\tfrac{1}{1-z_1^2}$, $z_k$, $\tilde{\vartheta}_i$, $\hat{\vartheta}_i$, $\omega_j$, $\tilde{\gamma}_j$, and $ \hat{\gamma}_j, (k=2, \ldots, n, i=1, \ldots, n, j=1, \ldots , n-1)$ are bounded. \textcolor{black}{Since $|y(0) - y_d(0)| < \rho(0)$,} one has that $z_1$ always evolves within the set $\Omega_z= \{z_1~|~|z_1|<1\}$. Then, $x_i$, $r$, $\lambda$, $\zeta_i$ and $\varphi_i$ are bounded. \textcolor{black}{From Remark \ref{remark:xi}, $\xi$ is bounded.} Moreover, since the control gains $\sigma_1$ and $\sigma_2$ are bounded based on their definition, then the virtual control inputs $\alpha_j~(j=1, 2, \ldots , n-1)$ and the actual control input $u$ remain bounded during the entire process of system operation. Since $\omega_j, \alpha_j \in \mathcal{L}_\infty$, then $\alpha_j^c$, $\dot{\alpha}_j^c$, the parameters estimate $\Hat{\vartheta}_i, \hat{\gamma}_j$  are bounded. Therefore, all signals in the closed-loop system remain bounded. 2) Since $\log\tfrac{1}{1-z_1^2}$, then $|z_1|<1$ and the output always evolves within the following prescribed set $\Omega_y=\{(e, \rho)~|~|e(t)|<\rho(t)\}$. 3) To prove practical prescribed-time stability and asymptotic convergence, the following two cases are considered.
        \vspace{-0.1in}
        
\textcolor{black}{\textit{Case 1:} For $0 \leq t < T$, we have
\begin{align}\label{eq:Vndot1}
    \dot{V}_n \leq - \chi_1 V_n + \chi_2,
\end{align}
where $\chi_1 = \bar{\varsigma} \sigma_1(t)$, $\bar{\varsigma} = 2 \min\{\varsigma_{z_i}, \varsigma_{\omega_i}, \iota_{\vartheta_i}, \iota_{\gamma_i}\}$, and $\chi_2 = \sum_{i=1}^n \sigma_2 \vartheta_i^2 + \sum_{i=1}^{n-1} \sigma_2 \gamma_i^2 + \tilde{D}$ such that $|\chi_2| \leq \bar{\chi} < \infty$. Let $\eta > 0$ be a constant such that $V_n(0) \leq \eta$. From \eqref{eq:Vndot1}, it follows that $\dot{V}_n < 0$ on $V_n = \eta$ when $\chi_1 > \frac{\chi_2}{\eta}$, which means that $V_n(t) \leq \eta$ for all $0 \leq t < T$ and $V_n(0) \leq \eta$.} \textcolor{black}{Moreover, the residual set of the solution of the system is $\lim_{t\to T} V_n(t)\leq \Omega$  where $\Omega = \exp(-\bar{\varsigma}\bar{\sigma} T)V_n(0) + \tfrac{\bar{\chi}}{\bar{\sigma}}(1-\exp(-\bar{\varsigma}\bar{\sigma} T))$. From the definition of $V_n(t)$, we have $|z_1|<\min\{\sqrt{1-10^{-2\underline{g}_1\Omega}}, \rho_T\}$ and $|z_i|<\sqrt{2\underline{g}_i\Omega}$ for $i=2, 3, \ldots , n$. Thus, a larger \( \bar{\sigma} \) value results in a smaller convergence region, which improves control precision.}

\textit{Case 2:} \textcolor{black}{For $t \geq T$, we have
  $\sigma_1(t) = \bar{\sigma}$  and
  $\sigma_2(t) = \bar{\sigma} \varepsilon(t-T)$. Since $\lambda=\frac{1}{1-z_1^2}$, from Lemma \ref{lem:log}, we obtain that $\log \frac{1}{1-z_1^2} <\frac{z_1^2}{1-z_1^2} =\lambda z_1^2 $ in the set $|z_1|<1$ and $-\frac{\varsigma_{z_1} \bar{\sigma}}{\underline{g}_1} \lambda z_1^2 < -\frac{\varsigma_{z_1} \bar{\sigma}}{\underline{g}_1} \log\left(\frac{1}{1 - z_1^2}\right).$ Moreover, since $- \log\frac{1}{1 - z_1^2} \leq - z_1^2$ holds, then the inequality \eqref{eq:Vndot2} is rewritten as}
\begin{align} \label{eq:Vndot3}
    \textcolor{black}{\dot{V}_n \leq} & \textcolor{black}{-\sum_{i=1}^{n}\tfrac{\varsigma_{z_i}\bar\sigma}{\underline{g}_i} z_i^2 -\sum_{i=1}^{n-1} \varsigma_{\omega_i}\bar\sigma \omega_i^2  + \tilde{D}} \nonumber \\
& \textcolor{black}{-\bar{\sigma} \varepsilon(t-T)\left(\sum_{i=1}^{n} (\tilde{\vartheta}_i^2-{\vartheta}_i^2)+ \sum_{i=1}^{n-1}  (\tilde{\gamma}_i^2- {\gamma}_i^2) \right) .} 
\end{align}
Integrating the inequality \eqref{eq:Vndot3} from $T$ to $t$ gives
\begin{align}\label{eq:aaa}
    V_n(t) &\leq  V_n(T) - \bar{\sigma}\int_{T}^t \left( \sum_{i=1}^n \tfrac{\varsigma_{z_i}}{\underline{g}_i} z_i^2(s) \, ds + \sum_{i=1}^{n-1} \varsigma_{\omega_i} \omega_i^2(s) \right) \nonumber \\
    &+ \bar{\sigma} \int_{T}^t\left( \sum_{i=1}^n \varepsilon(s-T) \vartheta_i^2 +\sum_{i=1}^{n-1}  \varepsilon(s-T) \gamma_i^2\right)ds \nonumber \\
    &+ \int_{T}^t \tilde{D}(t_0, s) \, ds.
\end{align}
Then, \eqref{eq:aaa} can be further represented as
\begin{align*}
\lim_{t\to\infty}&\bar{\sigma}\int_{T}^t \Big(\sum_{i=1}^n \tfrac{\varsigma_{z_i}}{\underline{g}_i} z_i^2(s) + \sum_{i=1}^{n-1} \varsigma_{\omega_i} \omega_i^2(s)\Big)ds \leq V_n(T) \nonumber \\ 
    &  +\lim_{t\to\infty} \bar{\sigma} \int_{T}^t\Big( \sum_{i=1}^n \varepsilon(s-T) \vartheta_i^2 +\sum_{i=1}^{n-1} \varepsilon(s-T) \gamma_i^2\Big)ds  \nonumber\\
    & +\lim_{t\to\infty} \int_{T}^t \tilde{D}(t_0, s) ds<+\infty.
\end{align*}
Note that $\sum_{i=1}^{n}(p \psi_{i2} \circ \varpi_1^{-1}(2D))^2=0$ after $t=T^0$. Since $z_i(t),~\omega_i(t)\in\mathcal{L}_\infty \cap \mathcal{L}_2$ and $\dot{z}_i(t),~\dot{\omega}_i(t)\in\mathcal{L}_\infty$, Barbalat's lemma implies that $\lim_{t\to\infty}\omega_i(t)=0$ and $\lim_{t\to\infty}z_i(t)=0$. Since $\rho(t)>0$, then $\lim_{t\to\infty}e(t)=0$. The proof is completed here.



    
    \begin{remark}    
        In contrast to \cite{luo2024practical}, which only focuses on achieving GTGP, the proposed approach not only provides the GTGP but also guarantees prescribed performance. While \cite{cao2022practical} also achieves the GTGP and prescribed performance, this method relies on the convergence time determined by the funnel boundary, resulting in a fast increase in control input at the settling time of the funnel boundary. This is due to the absence of a prescribed-time control strategy aimed at ensuring fast convergence of the states. Furthermore, our proposed approach, in comparison to both \cite{luo2024practical, cao2022practical}, accomplishes the GTGP and prescribed performance and guarantees asymptotic convergence even in the presence of nonvanishing dynamic uncertainties through the utilization of the new $\sigma$-modification control technique.
    \end{remark}


    \begin{remark}
        In the adaptive laws, we employed a novel $\sigma$-modification strategy, through introduction of the gain $\sigma_2(t)$ in the leakage term. This gain reaches $\Bar{\sigma}$ at the time $T$ to guarantee practical prescribed-time stability. If the leakage term was removed, large initial conditions would cause a rapid increase in $\Hat{\vartheta}_i$ and $\hat{{\gamma}}$, leading to a high amplitude on the control input \cite{li2016adaptive}. Further, since asymptotic stability is not recovered under standard $\sigma$-modification method due to the existence of leakage term, after the convergence time $T$, the leakage term gradually vanishes as time tends to infinity, ultimately achieving asymptotic convergence.
    \end{remark}

    \begin{remark}
        \textcolor{black}{The parameters $\varsigma_{z_i}$, \( \varsigma_{\omega_i} \), \( \bar{\sigma} \), and \( \varepsilon \) must be carefully tuned to balance tracking speed, control smoothness, and sensitivity to noise. A larger $\varsigma_{z_i}$, \( \varsigma_{\omega_i} \), and \( \bar{\sigma} \) result in faster tracking but also lead to larger control input and increased sensitivity to noise. On the other hand, if \( \varepsilon \) is too small, it results in a smaller tracking error but it can cause chattering due to discontinuous control actions. 
        Fundamentally, asymptotic convergence is achieved using a bounded rate function $\sigma$ and vanishing gains governed by \( \varepsilon \). However, very small values of \( \varepsilon \) can lead to chattering issues. To facilitate smoother implementation, several solutions can be considered, including carefully adjusting the convergence rate of \( \varepsilon \) or setting a very small lower bound for \( \varepsilon \) to ensure that the controller remains smooth at all times.}
    \end{remark}

    \begin{remark}
        \textcolor{black}{One application of prescribed-time prescribe performance control is in spacecraft attitude control for observation and communication missions, where precise orientation within a strict time interval is essential for pointing antennas, cameras, or sensors. For instance, in Earth observation satellites, any deviation in attitude can degrade image quality or cause data loss. Another important application is in robotic surgery, where surgical robots must perform precise movements within a defined time frame. In minimally invasive procedures, any deviation in precision or failure to meet prescribed performance bounds can jeopardize patient safety, making prescribed-time control with guaranteed performance crucial for successful outcomes.}
    \end{remark}

   \begin{remark}
 \textcolor{black}{The global finite-time result in \cite{wang2023adaptive} handles unmodeled dynamics by assuming they are exp-input-to-state practically finite-time stable, which is essential to break the circular dependence between $x_1$ and the dynamic signal $r$. Similarly, if the global prescribed-time scheme of \cite{shi2022prescribed} were extended to include unmodeled dynamics, one would need to assume that the unmodeled dynamics are exp-input-to-state practically prescribed-time stable. This is a much stronger and less realistic requirement than the standard exp-ISpS condition adopted in this paper. Although our result is semi-global, it relies only on the mild and widely used Assumption \ref{ass:exp}, reduces control complexity through a DSC framework, and ensures prescribed performance, bounded control, and asymptotic convergence in the presence of nonvanishing unmodeled dynamics. Furthermore, the semi-global region of attraction, defined by the bounds $\mathcal{B}_{1}$ and $q$, can be systematically enlarged by adjusting the parameters of the performance function $\rho(t)$ and the controller gains.}
\end{remark}


\section{Simulation results} \label{se:simulation}
We use the following two examples to evaluate the performance of the proposed control framework.

\begin{figure}
  \centering
  \includegraphics[width=0.8\columnwidth]{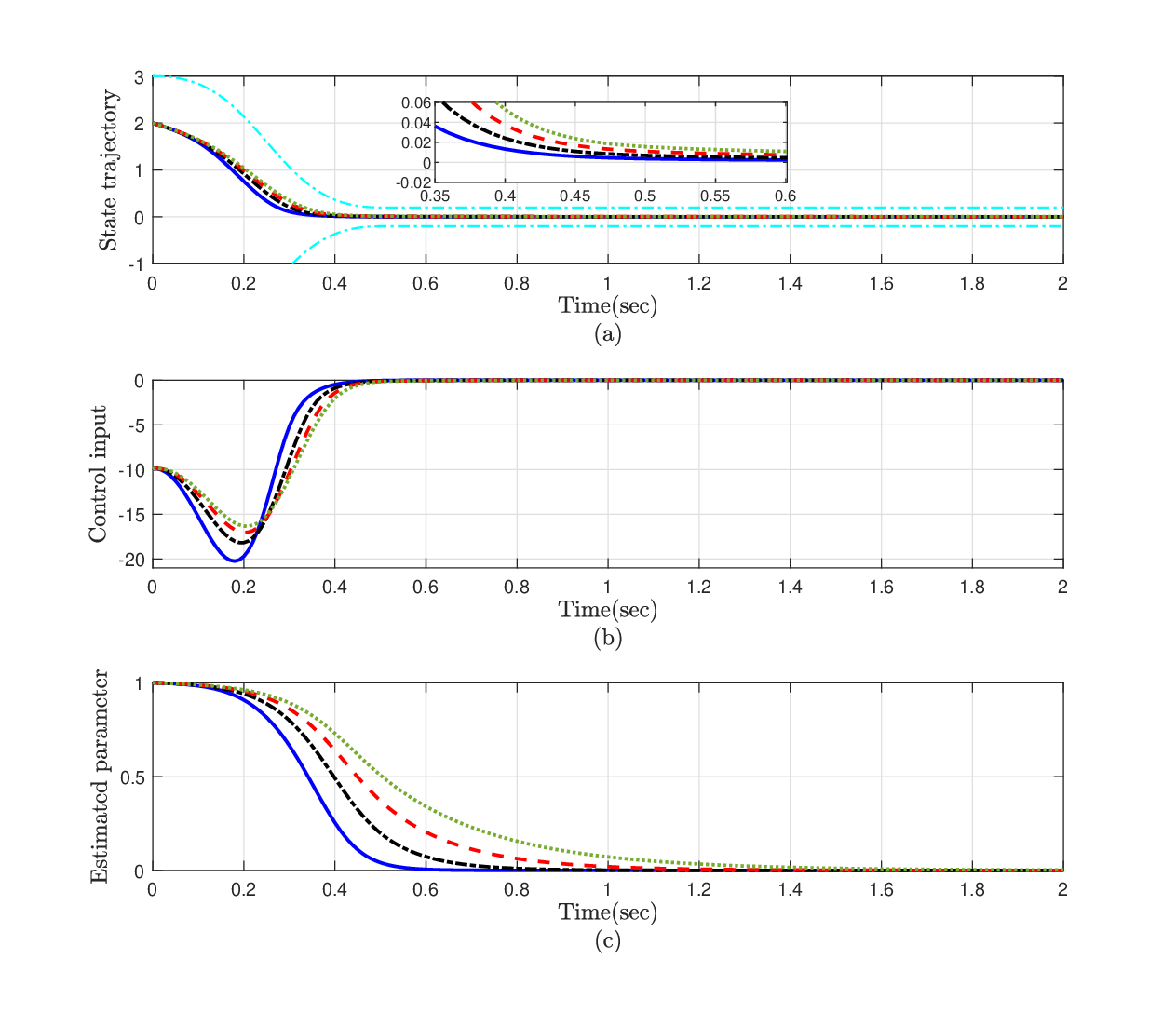}
  \caption{\textcolor{black}{(a) State trajectory; (b) control input; (c) estimated parameter ($\bar{\sigma}=100$: solid line, $\bar{\sigma}=50$: dash-dotted line, $\bar{\sigma}=30$: dashed line, $\bar{\sigma}=20$: dotted line)}} \label{fig:Ex1}
\end{figure}

\textbf{Example 1:} To show the effect of the rate function $\sigma(t)$, we consider the following uncertain first-order system:
\begin{align}\label{eq:first-order sys}
    \dot{x}(t)= g(t) u(t)+\theta(t) \phi(x),
\end{align}
where $g(t)=1-0.5\cos(tx),$ $\theta(t)= 1+0.5\sin(t)$, $\phi(x)=x \sin(x)$ and initial condition $x(0)=2$. The objective is to make the system state converge to zero in spite of the unknown control coefficient $g(t)$ and time-varying parameter $\theta(t)$. 
We apply the proposed prescribed-time controller to system \eqref{eq:first-order sys} for different values of $\bar{\sigma}=[20, 30, 50, 100]$. The rate function $\sigma(t)$ and the performance function $\rho(t)$ are constructed based on \eqref{eq:mu} with $T = 0.5$, $\upsilon_\sigma=0.4$, $\rho_0=3$, $\rho_T=0.2$, and $\upsilon_\rho=1$. The other control parameters are selected as $\varsigma_z =1$, $\varepsilon=\exp(-0.1t)$, and $\iota_\vartheta=0.1$.

From Fig. \ref{fig:Ex1}, the trajectory $x(t)$ evolves within the prescribed boundary and it converges to a small vicinity of the origin. A larger value of $\bar{\sigma}$ results in a smaller error at the prescribed time $T=0.5$. 
Since $\sigma(t)$ starts at $1$ for all cases, the initial control effort is the same for different values of $\bar{\sigma}$. \textcolor{black}{Moreover, the control signal is bounded for all $\bar{\sigma}$ and the method is valid for infinite-time interval. While larger $\bar{\sigma}$ results in a smaller tracking error at $t=T$, it requires more energy as shown in Table \ref{tab:energy}. Moreover, the update parameter remains bounded for all values of $\bar{\sigma}$.}

\begin{table}\caption{\textcolor{black}{Energy consumption for different $\bar{\sigma}$ in Example 1}}
    \centering 
    \begin{tabular}{c c c c c} \arrayrulecolor{red} \hline
        {} & \textcolor{black}{$\bar{\sigma}=20$} & \textcolor{black}{$\bar{\sigma}=30$} & \textcolor{black}{$\bar{\sigma}=50$} & \textcolor{black}{$\bar{\sigma}=100$} \\ \hline
        \textcolor{black}{Energy} & \textcolor{black}{$1066$} & \textcolor{black}{$1128$} & \textcolor{black}{$1231$} & \textcolor{black}{$1422$}\\ \hline
    \end{tabular}    
    \label{tab:energy}
\end{table}

\begin{figure}
  \centering
  \includegraphics[width=0.8\columnwidth]{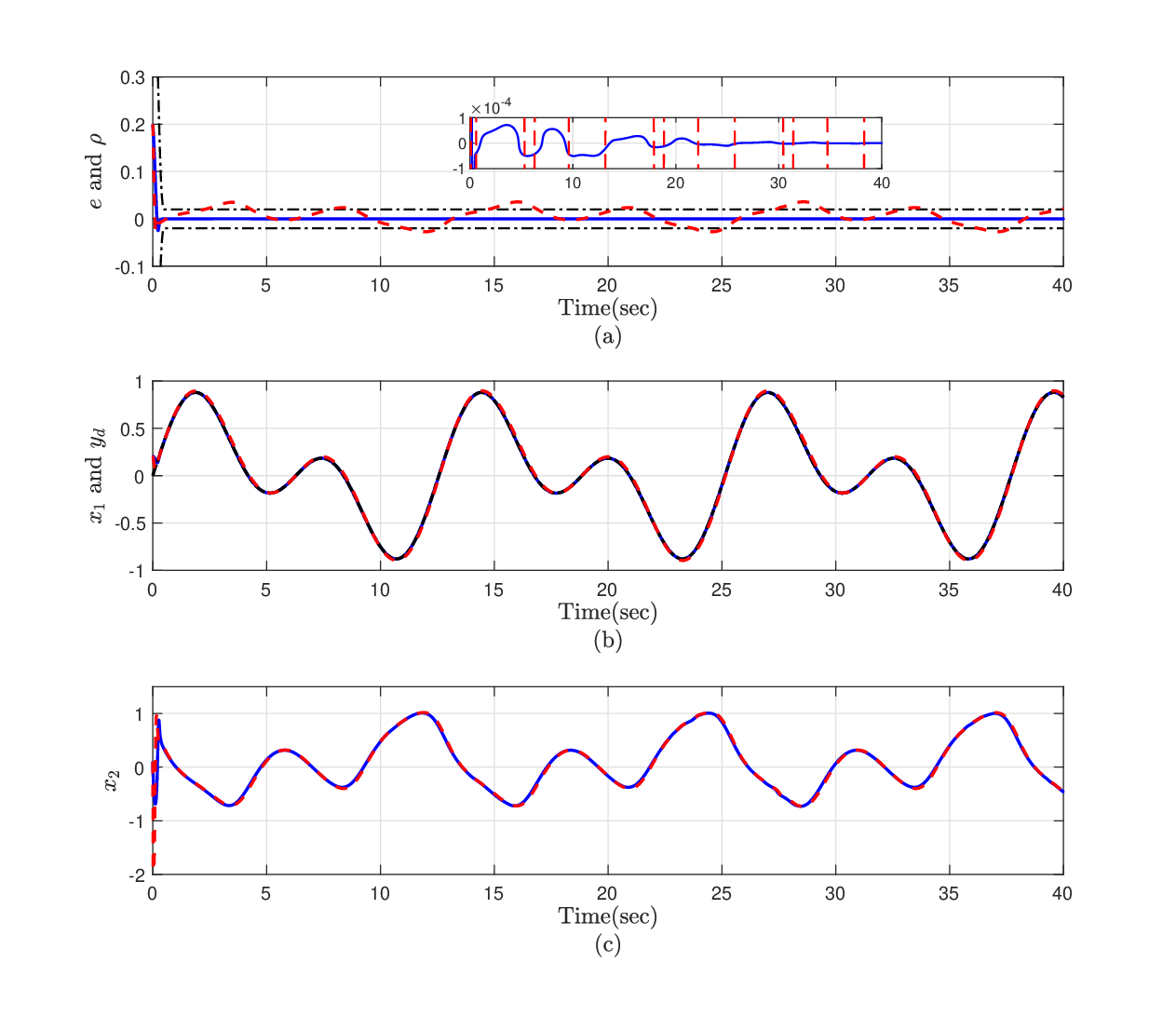}
  \caption{\textcolor{black}{(a) Output tracking error ($\rho(t)$: dash-dotted line, $e(t)$ under proposed: solid line, $e(t)$ under \cite{zhang2017adaptive}: dashed line); (b) First state ($y_d(t)$: dash-dotted line, $x_1(t)$ under proposed: solid line, $x_1(t)$ under \cite{zhang2017adaptive}: dashed line); (c) Second state ($x_2(t)$ under proposed: solid line, $x_2(t)$ under \cite{zhang2017adaptive}: dashed line)}}\label{fig:e_x1_x2}
\end{figure}

\begin{figure}
  \centering
  \includegraphics[width=0.8\columnwidth]{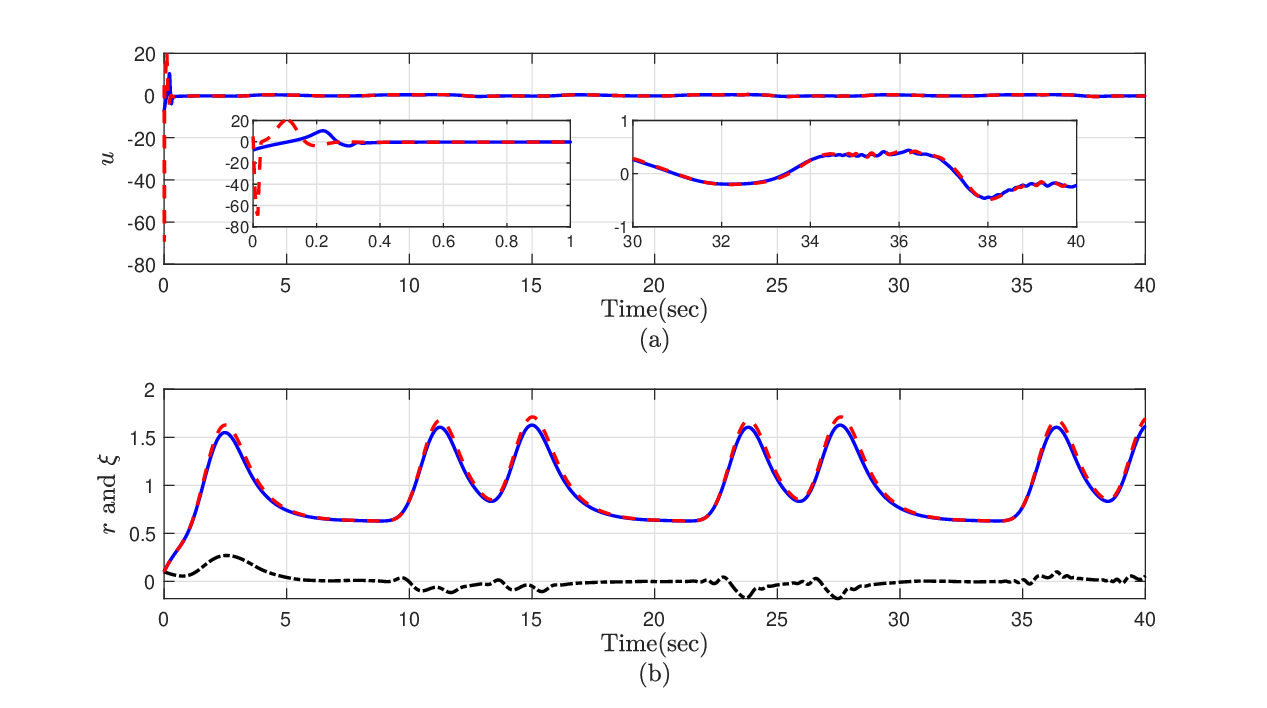}
  \caption{\textcolor{black}{(a) Control input ($u(t)$ under proposed: solid line, $u(t)$ under \cite{zhang2017adaptive}: dashed line); (b) dynamic signal ($\xi(t)$: dash-dotted line, $r(t)$ under proposed: solid line, $r(t)$ under \cite{zhang2017adaptive}: dashed line)}} \label{fig:u_r}
\end{figure}

\begin{figure}
  \centering
  \includegraphics[width=0.8\columnwidth]{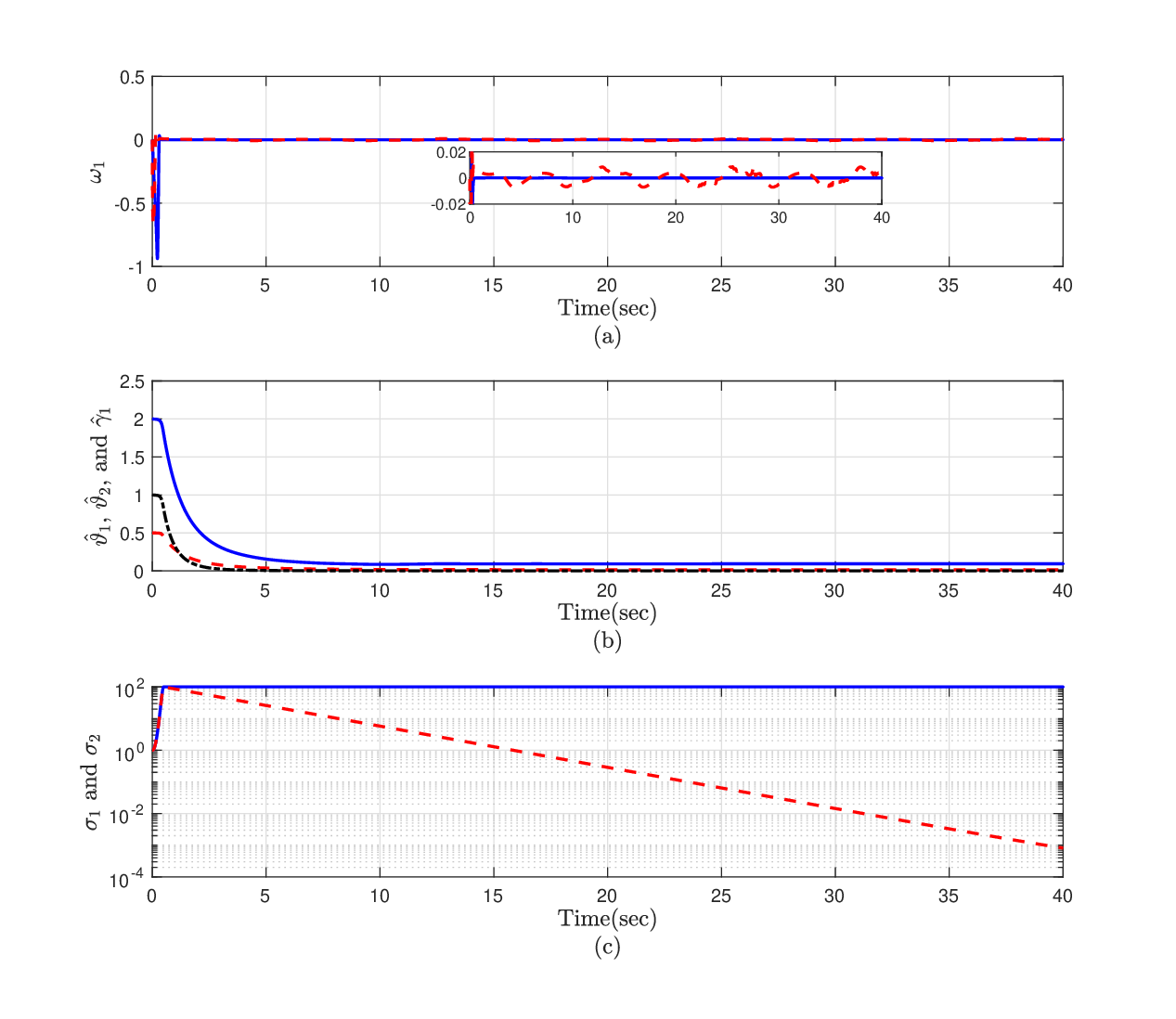}
  \caption{\textcolor{black}{(a) Filter error ($\omega_1(t)$ under proposed: solid line, $\omega_1(t)$ under \cite{zhang2017adaptive}: dashed line); (b) Estimated parameters ($\hat{\gamma}_1(t)$: dash-dotted line, $\hat{\vartheta}_1(t)$: solid line, $\hat{\vartheta}_2(t)$: dashed line); (c) Time-varying gains ($\sigma_1(t)$: solid line, $\sigma_2(t)$: dashed line)}} \label{fig:w_theta_sig}
\end{figure}

\textcolor{black}{\textbf{Example 2:} We consider the following nonlinear system with dynamic uncertainty (\cite{zhang2017adaptive})
$$\begin{cases}
        \dot{\xi} = -\xi+0.5x_1^2\sin(x_1t), & \\
        \dot{x}_1= x_1\exp(-0.5x_1)+(1+x_1^2)x_2+\delta_1(\xi,x_1,x_2,t), & \\
        \dot{x}_2= x_1x_2^2+(3-\cos(x_1x_2))u+\delta_2(\xi,x_1,x_2,t), & \\
        y = x_1
    \end{cases}$$
where $\delta_1(\xi,x_1,x_2,t) = 0.2\xi x_1\sin(x_2t)$ and $\delta_2(\xi,x_1,x_2,t)$ $=0.1\xi \cos(0.5x_2t)$ with $\psi_{11} = \sqrt{x_1^2+0.1}$, and $\psi_{12} = \psi_{22} = \sqrt{r^2+0.1}$. The desired trajectory and the dynamic signal are taken as $y_d(t)=0.5(\sin(t)+\sin(0.5t))$ and $\dot{r} = -r+2.5x_1^4+0.625$. The initial conditions of the states are set as $[\xi(0), x_1(0), x_2(0)]^\top = [0.1, 0.2, 0.1]^\top$. The proposed controller's parameters are taken as $T = 0.5$, $\bar{\sigma}=100$, $\varsigma_{z_1} = \varsigma_{z_2} =5$, $\varsigma_\omega =1$, $\varepsilon=\exp(-0.3t)$, $\upsilon_\sigma=0.2$, $\iota_\gamma=0.1$, $\iota_\vartheta=0.05$, $\rho_0=0.5$, $\rho_T=0.02$. We compared our approach with the adaptive DSC method for nonlinear systems with dynamic uncertainties and uncertain control gain, as described in \cite{zhang2017adaptive}. The controller parameters are considered the same as those in the original reference. The simulation results are shown in Figs. \ref{fig:e_x1_x2}--\ref{fig:w_theta_sig}.}

\textcolor{black}{From Fig. \ref{fig:e_x1_x2}. (a), the output tracking error under the proposed method converges to a prescribed region within the prescribed time $T$ and asymptotically converges to zero. However, Zhang’s method shows a persistent tracking error that does not vanish over time which results in less accurate performance. From Fig. \ref{fig:e_x1_x2}. (b), although both controllers exhibit good tracking performance, the proposed method achieves more precise control, as supported by the results in Fig. \ref{fig:e_x1_x2}. (a). From both Figs. \ref{fig:e_x1_x2}. (b) and (c), it is evident that the states remain bounded under both controllers, ensuring system stability.
Based on Fig. \ref{fig:u_r}. (a), the proposed controller requires less initial control effort due to the use of the rate function, which starts from 1 and gradually increases over time, ensuring given-time given-precision control. In contrast, Zhang’s controller requires a higher control gain to stabilize the tracking error which leads to a larger initial control effort. Fig. \ref{fig:u_r}. (b) shows that the unmodeled dynamics and dynamic signals remain bounded, and both controllers effectively handle the unmodeled dynamics.
From Fig. \ref{fig:w_theta_sig}. (a), the filter error under the proposed method is driven to a small region within the prescribed time and asymptotically converges to zero. However, it cannot vanish over time when Zhang's method is applied.
As it is observed from Fig. \ref{fig:w_theta_sig}. (b), the update parameters remain bounded, indicating stable adaptation of the system.
Fig.\ref{fig:w_theta_sig}. (c) illustrates the time-varying gains \( \sigma_1 \) and \( \sigma_2 \). Both start from 1 and reach their final value \( \bar{\sigma} \) within the prescribed time \( T \). After that, \( \sigma_1 \) remains constant, while \( \sigma_2 \) begins to vanish to provide asymptotic convergence of the errors.}

\section{Conclusion} \label{se:conclusion}
This paper proposes an adaptive practical prescribed-time control framework with prescribed performance for uncertain nonlinear systems subject to uncertain time-varying control coefficients as well as parametric and dynamic uncertainties. Specifically, we have developed an adaptive DSC scheme that effectively reduces the control complexity and ensures prescribed performance. By introducing a monotonically increasing rate function over $[0, T)$ that freezes at a user-specified time $T$ and a new ${\sigma}$-modification strategy in which the leakage term starts to vanish at $t\geq T$, our method guarantees prescribed-time convergence to a specific region and ultimately achieves asymptotic convergence. 
\textcolor{black}{In this paper, we assume that signs of the control coefficients \( g_i \) are known and all the states $x$ are available for feedback. Relaxing these assumptions will certainly increase applicability. Thus, working towards practical prescribed-time output feedback control with asymptotic convergence in the presence of totally unknown control coefficients deserves further investigation.}

\section*{Acknowledgment}
The authors would like to thank the reviewers and editors for their constructive suggestions that have significantly helped us improve the quality of this article.

\bibliographystyle{myagsm}        
\small{
\bibliography{Samplee}           
}

\end{document}